\documentclass{article}

\usepackage{arxiv}

\usepackage[utf8]{inputenc} 
\usepackage[T1]{fontenc}    
\usepackage{hyperref}       
\usepackage{url}            
\usepackage{booktabs}       
\usepackage{amsfonts}       
\usepackage{nicefrac}       
\usepackage{microtype}      
\usepackage{lipsum}
\usepackage{graphicx}
\graphicspath{ {./images/} }

\title{On the ubiquity of the ruler sequence}
\author{
 Juan Carlos Nu{\~{n}}o \\
  Department of Applied Mathematics\\ 
    Universidad Polit{\'e}cnica de Madrid\\ 
Madrid - 28040, Spain\\
  \texttt{juancarlos.nuno@upm.es} \\
   \And
 Francisco J. Mu{\~{n}}oz \\
\\
   \texttt{f.j.munoz.ortega@gmail.com}  \\
  \\
 \\
}

\begin{document}
\maketitle

\begin{abstract}
The {\it ruler function} or the {\it Gros sequence} is a classical infinite integer sequence that is underlying some interesting mathematical problems. In this paper, we provide four new problems containing this type of sequence: (i) a demographic discrete dynamical automata, (ii) the middle interval Cantor set, (iii) the construction by duplication of polygons and (iv) the horizontal visibility sequence at the accumulation point of the Feigenbaum cascade. In all of them, the infinte sequence is obtained by a recursive procedure of duplication. The properties of the ruler sequence, in particular, those relating to recursiveness and self-containing, are used to get a deeper understanding of these four problems. These new representations of the ruler sequence could inspire new studies in the field of discrete mathematics. 
\end{abstract}

\section{Introduction}

\begin{figure}
\centering
\includegraphics[width=0.9\textwidth]{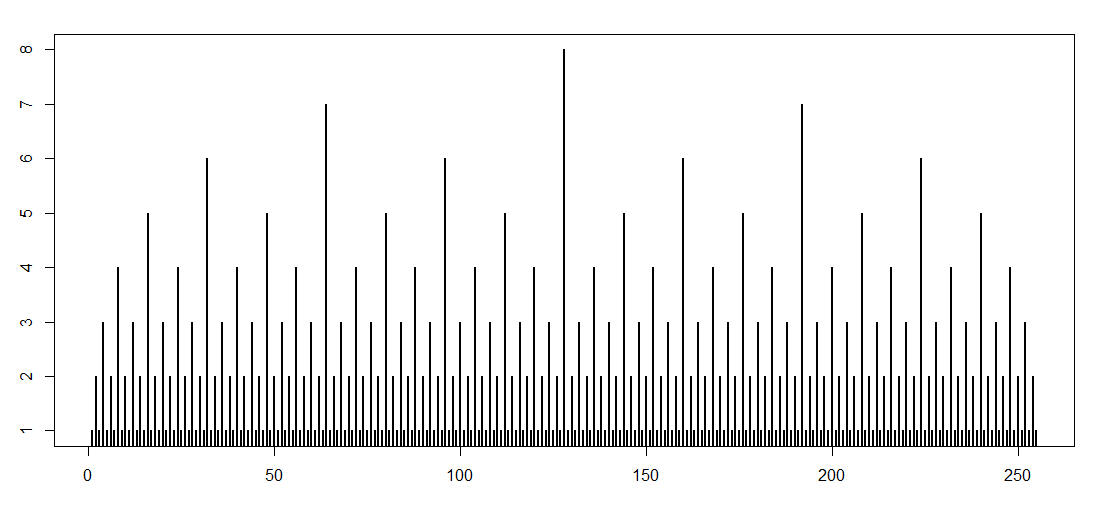}
\caption{Representation of the $2^8-1$ first terms of the ruler function. The X-axis represents the position in the sequence and the Y-axis the value of the term at each position. Note the form of the set of points that resemble the marks of a ruler. The points are generated from the algorithm presented in section \ref{someclass}, with $nmax=8$.}
\label{rulerfig}
\end{figure}

The ruler sequence, also known as the ruler function, is referred in N. J. A. Sloane, Sequence A001511 in OEIS (On-Line Encyclopedia of Integer Sequences) \cite{A001511}, where many of its definitions are presented.  This sequence is commonly defined as the infinite integer sequence whose $n$th-term, $a(n)$, is the highest power of 2 that divides $2\,n$, for $n=1,2,\ldots$.  It is also named as the {\it Gros sequence} by Hinz {\it et al.} \cite{Hanoi} because he handles the so called {\it Baguenodier} or Chinese Ring Puzzle problem (that seeks to know how to disentangle a set of rings linked together in a wire) and obtains, as a finite solution, the first terms of the ruler sequence. The first terms of this sequence are: $\{1, 2, 1, 3, 1, 2, 1, 4, 1, 2, 1, 3, 1, 2, 1,... \}$ and the $2^8 -1$ first terms of the sequence are represented in Fig. \ref{rulerfig}, which resembles the ticks on a ruler. In this paper, we present four additional problems across discrete mathematics where this sequence is obtained as a characteristic of their recursive nature:  (i) a demographic automata, (ii) the middle interval Cantor set, (iii) the construction by duplication of polygons and (iv) the horizontal visibility of the Feigenbaum cascade.

Many important sequences are obtained as solutions of discrete dynamical models \cite{Salinelli}, where the sequence represents the value of the variable, i.e. size of a population, at each time step. Given an initial condition, the sequence forms the trayectory of the initial value problem that, in principle, extends to infinite time steps. 
 Perhaps, the first historical example of such a discrete model is the one that gives rise to the Fibonacci sequence \cite{Roberts}, where the $n$th-term is obtained as sum of the two previous ones, i.e. $X(n) = X(n-1)+X(n-2)$ for $n>2$ and $X(1)=X(2)=1$.  Contrary to this simple model, many population models include the space variable and the age of the individuals\cite{Caswell}. A particular approach to handle with these spatial demographic models are the cellular automata \cite{Wolfram}

Sharing the same level of relevance with these discrete dynamical systems, there is the Cantor set in the framework of number theory, combinatorics and fractal theory \cite{Falconer, Peitgen}.  It is commonly said that the Cantor set contains all essential elements of these mathematical fields \cite{Vallin}. Starting from the unit interval, the classical ternary Cantor set is formed recursively by dividing at each step the resulting intervals into three subintervals and leaving out the central one, generating one of the most famous fractal sets, schemetatically shown in \ref{Fig1}. 

Polygons are one of the most studied planar figures. Behind the simple definiton of being bounded by three or more sides, they have fundamental geometric properties \cite{Barnes}. A polygon with $s$-vertices or sides will be also referred as to $s$-gon. The problem of polygon construction is an old one \cite{Boyer, Kuh}.The manipulation of polygons is also a popular issue in Geometry \cite{Conway}. Many problems seek to explore the polygons properties and to relate them with other geometric figures. For instance, {\it chopping polygons into triangles} and the related problem of frieze patterns, where collective properties of numbers appear from simple local rules \cite{Conway}. The duplication of sides of polygons inscribed within a circle of unit radius can be used to estimate, for instance, the value of $\pi$ \cite{Boyer}.  This procedure can be done using regular polygons, i.e. that have all sides and angles equal, but  this is not necessary.

As it has been stated above, the logistic map exhibits chaotic behaviour for certain values of the growth rate. In the classical model: $x_{n+1} = r \, x_n \, (1- x_n)$, as the growth rate $r$ increases from 1 to 4, chaos appears and disappears infinitely after the well studied bifurcations \cite{May, Strogatz, Peitgen}.  The bifurcation diagram has been reinterpreted by applying the horizontal visibility algorithm \cite{Luque}. Recently, we have discovered new patterns in the sequences derived from the horizontal visibility algorithm in the infinite cascades to chaos that appears in this bifurcation diagram \cite{Nuno20}.

What all four problems have in common is that they all give rise to the ruler or Gros sequence as a fingerprint of their properties. In the next section, the ruler sequence is derived from the basic principles that these four problems share. The paper is organized as follows:  a brief overview of the ruler function and its properties is presented in the next section. Section 3 presents an alternative procedure to obtain the ruler sequence by duplication. In the following sections, this procedure is applied to four problems where the ruler sequence underlies. Finally, section 8 contains some concluding remarks.

\section{Some classical descriptions of the ruler function}\label{someclass}

As mentioned in the introduction, the ruler sequence is known to be related to certain classical problems as the division of even numbers by powers of 2 or the {\it The Tower of Hanoi} \cite{Hanoi}. With regards to the first problem, the ruler sequence $a(n)$ is the $p$-adic valuation of  $2\,n$.  Concerning the second, the ruler sequence coincides with the sequence of movements that solves it,  if we numerate, in increasing order of diameter, the discs that are placed in one of the three pegs available, and codify the movements of the discs between the pegs in each step \cite{Hanoi}. According to \cite{Hanoi}, this sequence can be described recursively by de function $g$:

\[
g =
\bigg\{
\begin{array}{lcl}
1   &  if  &   k  \, odd \\
g_{k/2} + 1  &  if  &  k  \, even \\
\end{array}
\]
Another characterization of the ruler sequence reflecting the corresponding property of the Gray code, namely \cite{Hanoi}:
\[
\begin{array}{lclcl}
g_{2^n} & = & n +1  & if & n \in \mathbb{N}_0 \\
g_{2^n+k} & = & g_{2^n - k} & if & k \in [2^n - 1]
\end{array}
\]

A simple and useful way of generating the ruler sequence of size $nmax$ is given by this script \cite{Caroli}, written in R  \cite{R}:
\begin{verbatim}
nmax = 8;
r = 1;
for (n in 2:nmax){
  r =c(r,n,r)
  }
\end{verbatim}
Figure \ref{rulerfig} is drawn using this algorithm (with $nmax=8$).

Some authors treat  the ruler sequence and the Thomae's function as synonyms (see, for instance, \cite{Dunham}). Indeed, the ruler sequence can be viewed as a restriction of the Thomae's function to the dyadic rationals: those rational numbers whose denominators are powers of 2. If we define 
\[h(x) =
\bigg\{
\begin{array}{cl}
2^{-k}  & if   \hspace*{5mm} x= \frac{p}{2^k} \hspace*{5mm} with \hspace*{5mm} p \hspace*{5mm} odd   \\
0 & otherwise
\end{array}
\]
the ruler sequence coincides, after a linear translation, with the sequence of the exponents of the image of the set dyadic rationals by $h$, ordered from least to greatest.  For instance, if we restrict the dyadic rationals to those that have a power $n$ less than 4, the ordered image of $h$ is:
\[
\frac{1}{2^4},  \frac{1}{2^3}, \frac{1}{2^4}, \frac{1}{2^2}, \frac{1}{2^4}, \frac{1}{2^3}, \frac{1}{2^4}, \frac{1}{2}, \frac{1}{2^4}, \frac{1}{2^3}, \frac{1}{2^4}, \frac{1}{2^2}, \frac{1}{2^4}, \frac{1}{2^3}, \frac{1}{2^4}
\]
whose sequence of exponents is:
\[
-4, -3, -4, -2, -4, -3, -4, -1, -4, -3, -4, -2, -4, -3, -4
\]
This sequence can be moved to yield positive integers by adding the largest absolute number plus 1, in this case $4+1=5$, to all the members of the sequence:
\[
1, 2, 1, 3, 1, 2, 1, 4, 1, 2, 1, 3, 1, 2, 1
\]
The same transformation can be done for any $n$ and hence, obtain the infinite ruler sequence.

As referred in the introduction, another descriptions of the ruler function can be found in \cite{A001511}.  We end this section by pointing out two particular properties of the ruler sequence: (i)  it is self-contained \cite{Kimberling} and (ii) it is squarefree \cite{Allouche}. Both are consequences of the recursive properties of this sequence and of its fractal nature. These properties will be even more evident after taking into account the recursive definition and the problems we show in the next sections.

\section{An alternative derivation of the ruler sequence}\label{S1}

\begin{figure}
\centering
\includegraphics[width=0.9\textwidth]{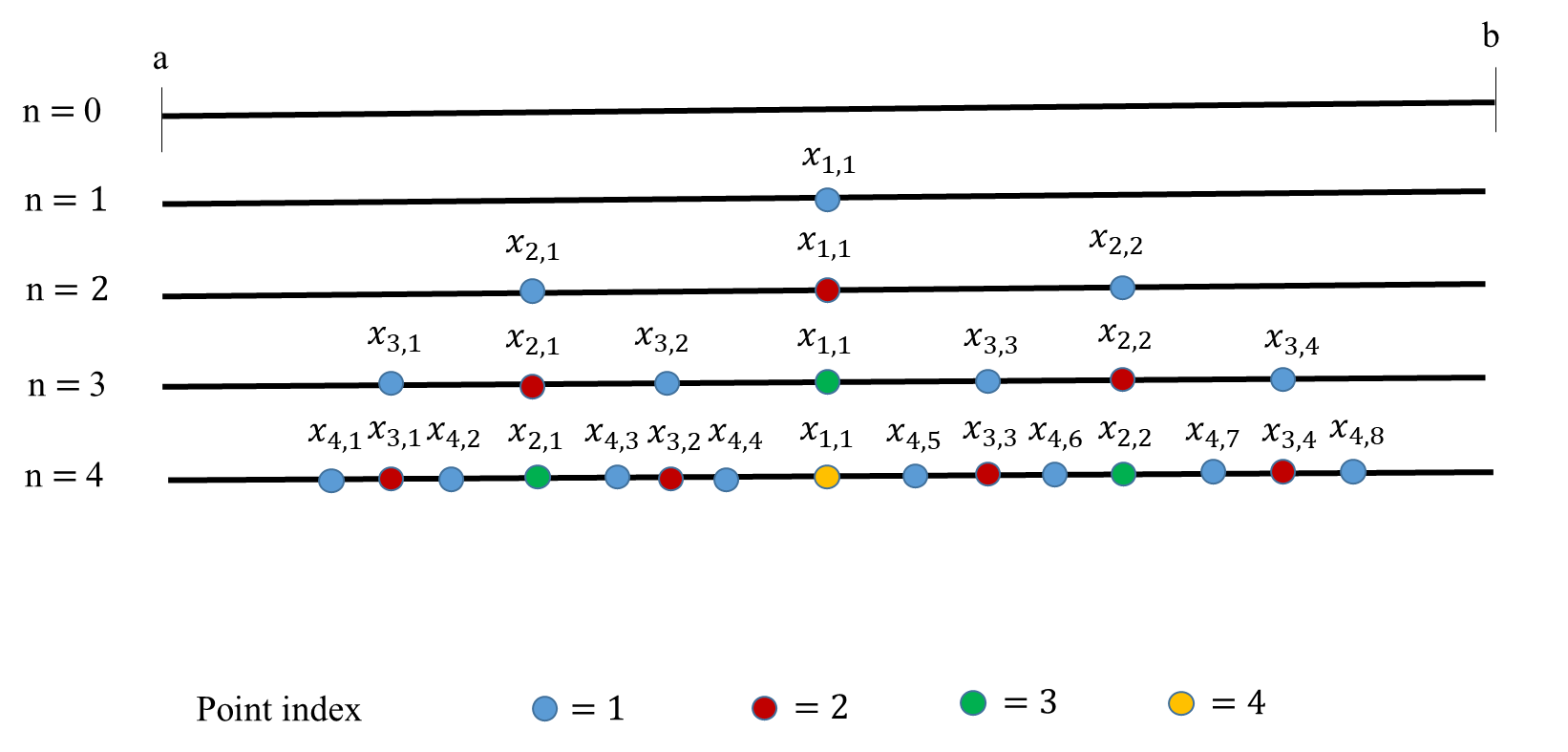}
\caption{Basic algorithm for the generation of the ruler sequence. From a given interval $I$, that can be the $(-\infty, \infty)$, a interior point $x_{1,1}$ is chosen and the following basic algorithm is applied. (i) insert two points at each side of $x_{1,1}$, (ii) for each point appeared in the previous step, $x_{n,k}$, insert two points in the partition at each side side in such a way that do not intersect with any other point that comes from other points of the partition. (iii) assign an index to each of the points of the new partition, 1 (blue) for the new points and increase in one unity for the rest of points. (iv) continue at infinitum.}
\label{FigB}
\end{figure}

Let us consider a sequence of successive partitions of a interval $I \in \mathbb{R}$,  finite or infinite, obtained by duplication of exclusively the points included in the previous step (see Fig. \ref{FigB}). If $x_{n,k}$ is the $k$th-point of the partition appeared in the step $n$, it gives rises to two new points in the next step $n+1$ that verify: 
\[
x_{n+1,i} < x_{n,k} < x_{n+1,j}
\]
In addition, both $x_{n+1,i}$ and $x_{n+1,j}$ must not overlap with the points inserted in the partition by duplication of other points, i.e. $x_{n+1,i-1} < x_{n+1,i}$ and $x_{n+1,j} < x_{n+1,j+1}$.  The other points already present in previous steps, but that do not duplicate, are also added to this next step $n+1$. 

The partition of the interval at each step is given by: 
\[
X_n = \{x_{n,1}, x_{n-1,1}, x_{n,2},  \ldots,x_{1,1} , \ldots, x_{n,2^{n-1}-1}, x_{n-1,2^{n-2}}, x_{n,2^{n-1}}\}
\]
For instance, as it can be seen in Fig.\ref{FigB}:
\[
X_1=\{x_{11}\}, \, X_2=\{x_{2,1},x_{11},x_{2,2}\}, \, X_3=\{x_{3,1},x_{2,1},x_{3,2},x_{1,1},x_{3,3},x_{2,2},x_{3,4}\}
\]
and 
\[
X_4=\{x_{4,1},x_{3,1},x_{4,2},x_{2,1},x_{4,3},x_{3,2},x_{4,4},x_{1,1},x_{4,5},x_{3,3},x_{4,6},x_{2,2},x_{4,7},x_{3,4},x_{4,8}\}
\]

The schematic representation of the sequence construction depicted in Fig. \ref{FigB} can be viewed as a tree, i.e. a graph without cycles. However, it had to be taken into account that, contrary to classical tree graphs, the nodes at the upper levels are also presented at lower levels with their corresponding indices. These trees has some similarities with the so called {\it Stern-Brocot trees}, that provides the positive rational numbers \cite{Graham}. 

The sequence of indices at step $n$, $r_n$, is formed by the index of each point of the partition at this step. Thus, at it can be seen in Fig. \ref{FigB}, $r_1 = \{1\}$, $r_2=\{1,2,1\}$ and $r_3=\{1,2,1,3,1,2,1\}$. Once the partition has been updated, we assign the corresponding index to all points of the partition generated at this step as follows: 1 for the new included points and increase by one unit for the rest of points that belong already to the partition. In so doing, the sequence of indices at step $n$ is:
\[
r_n = \{1, 2, 1, 3, \ldots, n , \ldots 3, 1, 2, 1 \}
\]
As $n$ tends to infinity, this sequence $r_n$ tends to the {\it ruler} sequence or Gros sequence \cite{A001511}.

It is a Regular Self-Containing sequence, i.e. it contains a proper subsequence that is identical to itself and, in addition, is regular \cite{Kimberling}. However, according to the definition given by Kimberling, it is not a fractal sequence \cite{Kimberling}. Note that this sequence verifies that after deleting the first occurrence of each positive integer the remaining sequence is the same as the original.

It can be shown that the asymptotic relative frequency of each natural number $k$ in the ruler sequence is $f(k) = 2^{-k}$. Obviously, the sum of these relative frequencies must be 1, i.e. 
\begin{equation}\label{agedistri}
\sum_k^{\infty} \frac{1}{2^k} = 1
\end{equation}

Let us now find the sum of the terms of the index sequence. For instance, $s_1=1$, $s_2= 4$ and $s_3=11$.  It is not difficult to find a recursive formula for this sequence:
\begin{equation}\label{rsum}
s_{n+1} = 2 \, s_n + n+1
\end{equation}
for $n=1,2,\ldots$. This expression can be rewritten in terms of $s_1=1$ as: 
\begin{equation}\label{explicitrsum}
s_{n+1} = 2^{n-1} + n+1 + 2 \, n + 2^2 \, (n-1) + \ldots +2^{n-2} \, 2
\end{equation}
for $n=1,2,\ldots$. 

As the partition is constructed, the size of the succesive terms of the index sequence increases.  It is straightforward to prove that  the number of terms in the sequence at step $n$ is:
\begin{equation}\label{expgrowth}
\#r(n) \equiv N(n) =  2^n -1 
\end{equation}
which provide an exponential growth.


It is also evident that the maximum integer in the $n$th-term is $max(r_n) = n$.  The ratio between the maximun integer and the length of the $n$th-term
\[
\lambda(n) = \frac{max(r_n)}{\#r(n)}
\]
can be expressed as a function of $N$ as:
\[
\lambda(N) = \frac{\log_2(N+1)}{N}
\]
This ratio behaves as $N^{-1}$ as $N$ tends to infinity.



\section{The solution of a recursive demographic automata}\label{automata}

The procedure of point generation presented in the previous section resembles the way individuals of a population self-reproduce. The points could be considered as individuals that are located in a one-dimensional array. Initially, a unique individual located in a spatial position gives birth to two new individuals that are placed at each side of him. According to the rules of duplication, individuals replicate only once,  just in the next step after they are formed.  Besides, the model considers that these indivuals are immortal, i.e. no rate death exists. Then, at each time step, this population located in the one-dimensional array increases its size  and the age of the individuals of the population at step $n$ is given by the sequence index $r_n$. That is, the ruler or Gros sequence would represent the age distribution of the population located in the one-dimensional array at infinite time (and infinite space).

The total population at step $n$ is given by equation \ref{expgrowth}, which represents an exponential growth. This dynamics can be obtained as the solution of the discrete equation:
\begin{equation}\label{discreteq}
N(n+2) = 2 \, \left(N(n+1) - N(n) \right) + N(n)  
\end{equation}
for $n=1,\ldots$ and 
\[
N(1) = 1; \, N(2) = 2 \, N(1) + 1= 3
\]
The first term of equation \ref{discreteq} represents the duplication of the individuals born in the previous step. The second term takes into account that all individuals of the previous step survive in the next, with an age increased by one. 

Equation \ref{discreteq} can also be rewritten as:
\begin{equation}\label{discreteq2}
N(n+2) = 3 \, N(n+1) - 2 \, N(n) 
\end{equation}
taking the form of the, so called, Horadam sequences that includes, among others, the Fibonacci and Lucas sequences \cite{Larcombe}. 

It is interesting to note that these equations are equivalent to the discrete equation:
\begin{equation}
N(n+1) = 2 \, N(n) +1 
\end{equation}
with the initial term $N(1)=1$, that represents the number of movements required to move completely $n$ discs from one peg to other, according to the rules of the towers of Hanoi game \cite{Hanoi}.

It is not difficult to prove that the difference $N(n+1) - N(n) = 2^{n}$, which means that, at each step $n$, $2^{n}$ individuals of the population are able to have offsprings, while the rest of the population, $2^{n} - 1$ individuals, increase their age by one unity. 

\begin{figure}
\centering
\includegraphics[width=0.6\textwidth]{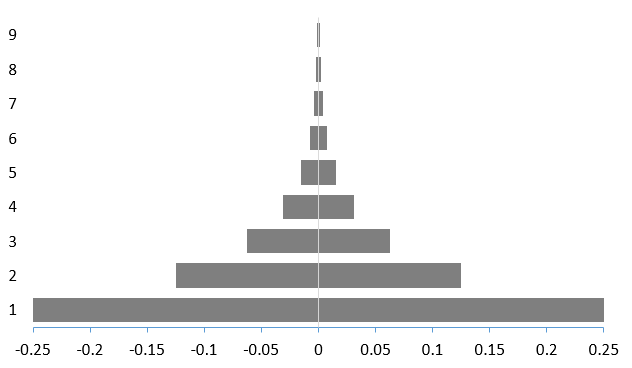}
\caption{Demographic pyramid obtained from the cellular automata model.  The lowest level is formed by newborn individuals, label by 1. The second level is formed by individuals of age 2 and, it can be seen, its proportion is half the proportion 1 age population. In general, the proportion of individuals of age $k$ is given by $p(k)=2^{-k}$.}
\label{pyr}
\end{figure}
\begin{figure}
\centering
\includegraphics[width=1.0\textwidth]{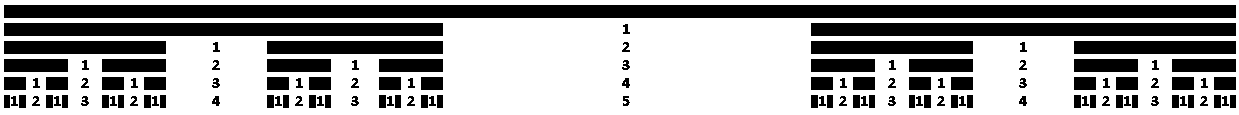}
\caption{Recursive generation of the Cantor set and  the corresponding index sequence. As the Cantor set is generated, the index sequence $r_n$ is obtained by recurrence:  after the first partition of the unit interval, for $n =1$, the index of the interval $i_{11}=(\frac{1}{3}, \frac{2}{3})$ is $i_{11}=1$ and the index sequence is $r_1 =\{1\}$.  Two new intervals appear in the next step, $n=2$, $i_{21} =(\frac{1}{9},\frac{2}{9})$ and $i_{22} =(\frac{7}{9},\frac{8}{9})$, whose indices are $i_{21}=i_{23}=1$. The index of the middle interval is increased by 1, $i_{11} =2$. Thus, at this step, the set sequence is: $C_2=\{i_{21},i_{11},i_{22}\}$ and the corresponding index sequence: $r_2=\{1,2,1\}$. As it can be seen, tis procedure continues in the next steps, forming the succesive terms of the Cantor set and its associated index terms.  The fractal structure of the set is evident.  The index sequence of the Cantor set is described by an array of an infinite number of columns of infinite length which contains the set of natural numbers. It is easy to see that the sequence is self-containing by construction: if we start from any 1, in any of the columns, the same sequence is obtained.}
\label{Fig1}
\end{figure}

The age distribution of the population is represented by expression \ref{agedistri} and it is depicted in the demographic pyramid (Fig. \ref{pyr}).

The order of the sequences $r_n$, as they are generated step by step, is consequence of the spatial disposition defined  in the cellular automata (see,  for instance, \cite{Cobeli}). Since the population growth is not bounded then, an array with infinite components is also needed. 

The demographic model can be modified to consider the death of individuals at a given age. For example, it could be thought that individual lifespan is divided into three periods: (i) fertility, (ii) maturity and (iii) senescence. Thus, after three time steps, individuals die and dissapear from the population. Accordingly, the discrete model \ref{discreteq} modifies to:
\begin{equation}\label{withdeath}
N(n) =
\bigg\{
\begin{array}{clc}
2^n - 1 & if & n=1,2 \\
7 \, 2^{n-3} & if & \hspace*{8mm} n =3, 4, 5, \ldots
\end{array}  
\end{equation}
The expression of equation \ref{withdeath} for $n=3, 4, \ldots$ is the solution of the geometric equation:
\[
N(n+1) = 2 N (n)
\]
 which means that, for this particular death rate, the population duplicates at each step. It is curious to recognize that, even when all the individuals finally die, the population size grows in an unlimited way. Therefore, to have a bounded population it is necessary to include an additional assumption that, either reduces the growth rate or increases the rate of disappearance of individuals. Obviously, in any case, the age distribution is not longer described by the ruler sequence.

\section{Middle interval indices of the Cantor set}

The Cantor set is a fractal subset of the interval $[0,1]$ \cite{Vallin, Falconer, Peitgen}. Geometrically, starting from the unit interval,  the classical Cantor set
obtains three equal subintervals of length $l_1=\frac{1}{3}$: two at each side of the middle open interval: $i_{11}= (\frac{1}{3}, \frac{2}{3})$, respectively: $[0, \frac{1}{3}]$ and $[\frac{2}{3}, 1]$.  We define the first term of the sequence of middle intervals as $C_1=\{i_{11}\}$ and the index of $i_{11}$ as $r_{1} = 1$, meaning  the number of steps since $i_{11}$ has been formed.  At the next step, only the first and the third intervals are partitioned into three new subintervals. The previous middle interval remains undivided forever. Thus, at this step, the unit interval is divided into 6 new subintervals of length $l_2=\frac{1}{9}$, whereas the previous middle interval has triple length:
$[0, \frac{1}{9}],  (\frac{1}{9}, \frac{2}{9}),  [\frac{2}{9}, \frac{1}{3}], (\frac{1}{3}, \frac{2}{3}), [\frac{2}{3}, \frac{7}{9}], (\frac{7}{9}, \frac{8}{9}), [\frac{8}{9}, 1]$.
The three middle intervals are: $i_{21} =  (\frac{1}{9}, \frac{2}{9})$,  $i_{11}= (\frac{1}{3}, \frac{2}{3})$ and $i_{22}= (\frac{7}{9}, \frac{8}{9})$. The indices of these three intervals are: $i_{21}= i_{22} = 1$, meaning that it is the first step they appear and $i_{11} = 2$, since it exists after 2 steps. The sequence of middle intervals at this step is given by $C_2=\{i_{21}, i_{11}, i_{22}\}$ and its corresponding sequence index by $r_2= \{1,2,1\}$.

At the third step, four new middle open intervals appear: $i_{31} = (\frac{1}{3^3},\frac{2}{3^3})$,  $i_{32} = (\frac{7}{3^3},\frac{8}{3^3})$, $i_{33} = (\frac{17}{3^3},\frac{18}{3^3})$ and $i_{34} = (\frac{25}{3^3},\frac{26}{3^3})$, that,  together with the previous $i_{21}, i_{11}$ and $i_{22}$, yields the 3th-term of the sequence of middle intervals:
\[
C_3 = \{i_{31},i_{21},i_{32}, i_{11}, i_{33}, i_{22}, i_{34}\}
\]
and the associated index sequence:
\[
r_3 =\{1,2,1,3,1,2,1\}
\]

For each step, the $n$th-term of the sequence of middle intervals, $C_n$, has $2^{n} - 1$ elements ($n=1,2,\ldots$) and so has the respective index sequence $r_n$. A general formula for the $n$th-terms of the sequence of intervals and the corresponding index sequence can be written as:
\[
C_n = \{i_{n1}, i_{(n-1) 1}, i_{n2}, i_{(n-2) 1}, \ldots, i_{11} , \ldots, i_{(n-1)2^{(n-2)}}, i_{n2^{n-1}}\}
\]
and the respective index sequence:
\[
r_n = \{1, 2, 1, 3, \ldots, n , \ldots 3, 1, 2, 1 \}
\]

Note that the construction of the sequences $r_n$ does not depend on the way the intervals are chosen. It only depends on the number of divisions that, as in this case, must be three for step, with a middle interval that remains undivided. Therefore, as stated in section \ref{S1}, the ruler sequence $r$ appears exclusively by the process of subdivision according to the basic rules there defined.

The distribution of interval lengths is provided by the respective index sequence, $r_n$. For example, in the previous case, $n=3$, there are 4 intervals of length $\frac{1}{3^3}$, 2 of length $\frac{1}{3^2}$ and 1 of length $\frac{1}{3}$. Thus, it is not difficult to obtain the whole length of the middle intervals that form the $n$th-term of the sequence of intervals:
\[
L_n = \sum_k \, l_n \, 3^{r_n(k)-1}
\]
where $l_n$ is the length of the $n$th-subdivision, i.e. $l_n=3^{-n}$ and the exponent $r_n(k)$ must be interpreted as the $k$th coordinate of $r_n$. Hence, this expression can be simplified to:
\[
L_n = \sum_k 3^{r_n(k)-(n+1)} 
\]
 For instance, the first 3 whole middle interval lengths are:
\[
L_1 = \frac{1}{3} ; L_2 = \frac{1}{3^{-2}}+ \frac{1}{3^{-1}}+ \frac{1}{3^{-2} }= \frac{5}{9}
\]
\[
L_3 = \frac{1}{3^{-3+0}}+ \frac{1}{3^{-3+1}}+ \frac{1}{3^{-3+0}} + \frac{1}{3^{-3+3}} + \frac{1}{3^{-3+0}} + \frac{1}{3^{-3+1}} + \frac{1}{3^{-3+0}}  = \frac{19}{27}
\]
This last expression can be also rewritten as:
\[
L_3 = 4 \frac{1}{3^3} + 2 \frac{1}{3^2} + \frac{1}{3}
\]

It is known that this sequence of the whole length tends to 1 as $n$ increases \cite{Vallin}.

\section{Vertex indices of infinite sided polygons}

\begin{figure}
\centering
\includegraphics[width=0.8\textwidth]{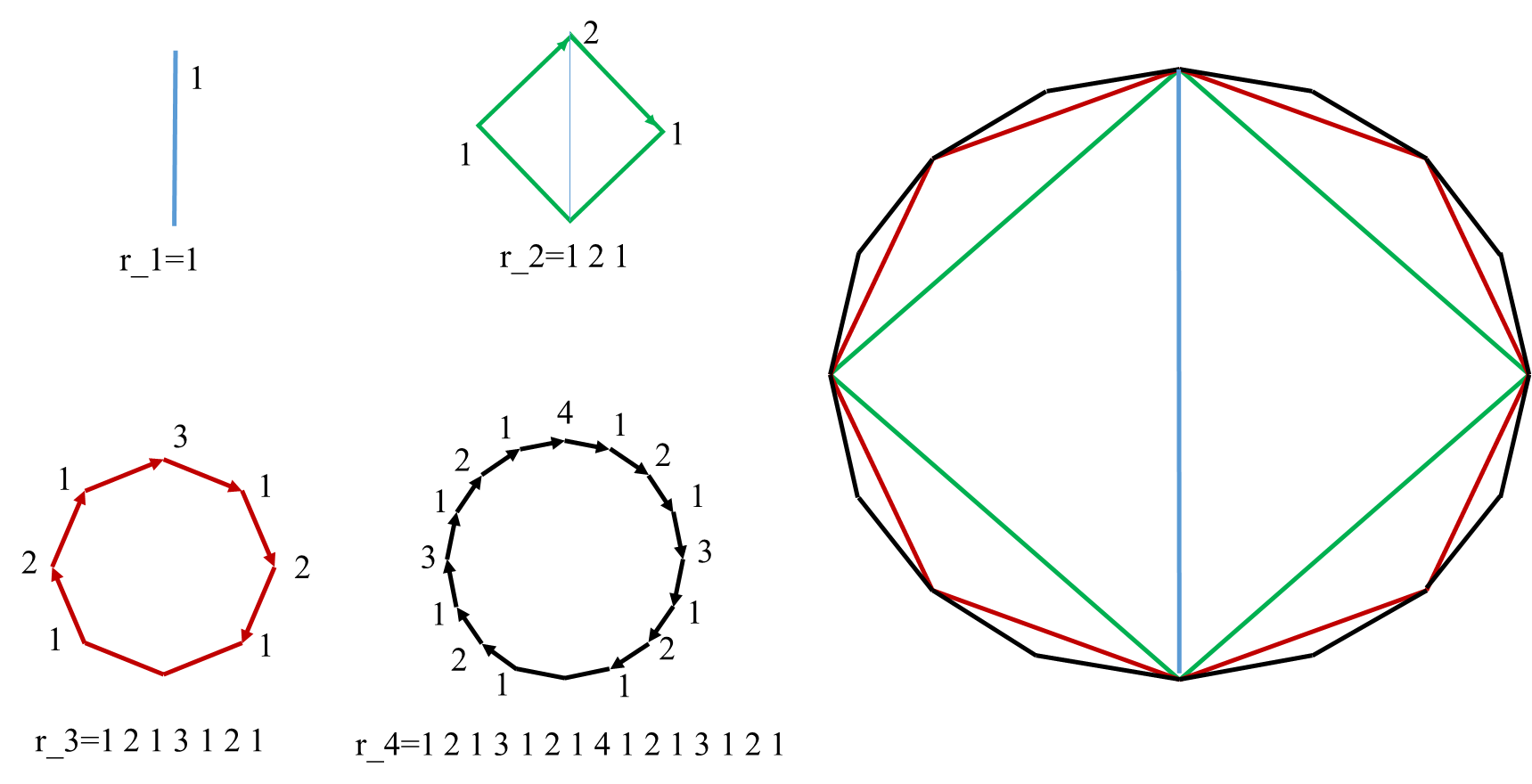}
\caption{Sequence generation from the regular polygons with $s=2^n$ vertices for $n=1,2,3$ and $4$. The next polygons with the double of sides (and vertices) as the previous one is constructed via bisecting the angles. The vertex index, that coincides with the number polygons to which this point belongs, is depicted for each polygon. }
\label{ngon}
\end{figure}

As stated in the introduction, the construction of polygons of any number of sides is a classical problem that has given rise to fundamentals results in geometry and number theory \cite{Boyer, Kuh}. For the sake of simplicity, let us consider the regular $s$-gons of $s$ vertices: $s=m \, 2^{n-1}$ for $n=1,2,…$  and $m$ a prime number. For instance,  if we divide equally each of the sides,  we obtain the sequence of $s$-gons with $s=2^n$ vertices  (see Figure \ref{ngon}). Note that for each side division, the same number of vertices appear. 

For each $s$, we define the {\it vertex index} as the number of $p$-gons with $p \leq s$ that share this vertex. For instance, as it can be seen in \ref{ngon}, the nothernmost vertex has index 4 because it belong to all $s$-gons for $s=1,2,3,4$.  At each step, the new generated vertices have always index 1.  For each $s$-gon, the sequence of vertex index $r_s$ is formed from $s-1$ vertices (the southermost vertex is not included), taken in a determined order, concretely from vertex located at the left of the southermost vertex and clockwise  (see Fig. \ref{ngon}).  For example, the vertex index sequence for $s=3$ then, $r_3=\{1, 2, 1, 3, 1, 2, 1, 4, 1, 2, 1, 3, 1, 2, 1\}$, which coincides with the same term of the ruler sequence. By construction, when $s$ tends to infinity, $r_s$ tends to the infinite ruler sequence. 

It is worthy to remark the recursive procedure that is behind the generation of the ruler sequence: each $s$-gon inscribes all the  $p$-gons with $p < s$. This property is reflected in the recursive pattern of the ruler function infinite sequence. Moreover, as in the previous cases, the process of duplication of sides must not be necessarily regular, i.e. the length of the sides has not to be equal.  Equivalently, this means that the position of the vertices on the circle has not to be equally spaced but, as stated in section \ref{S1}, they never must cross each other.

\section{Forward horizontal visibility at the accumulation point of the Feigenbaum cascade}

\begin{figure}
\centering
\includegraphics[width=0.8\textwidth]{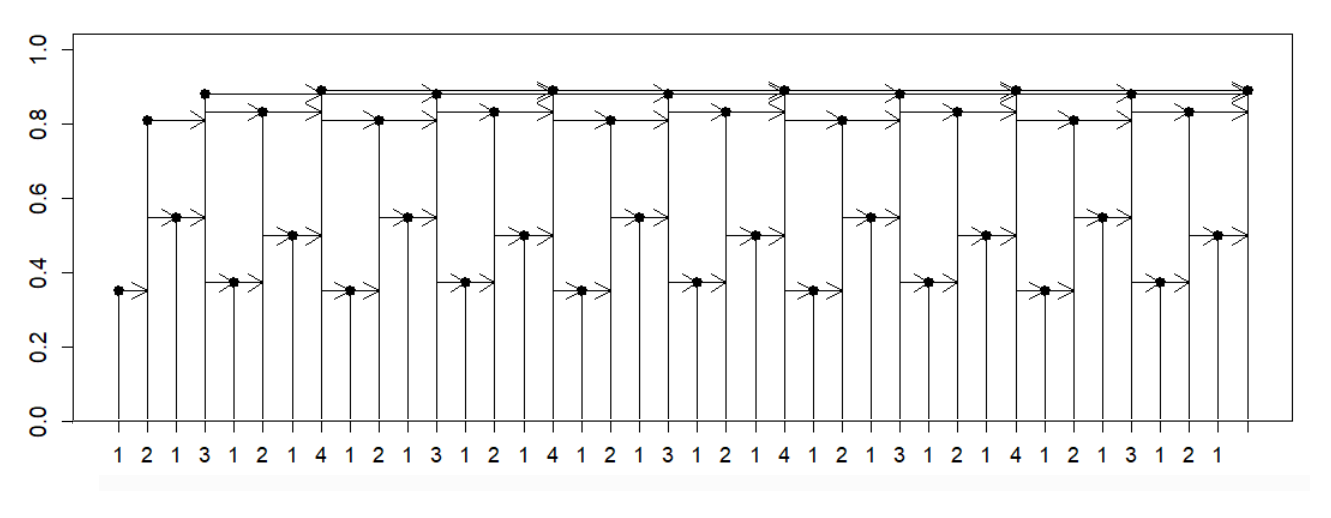}
\includegraphics[width=0.8\textwidth]{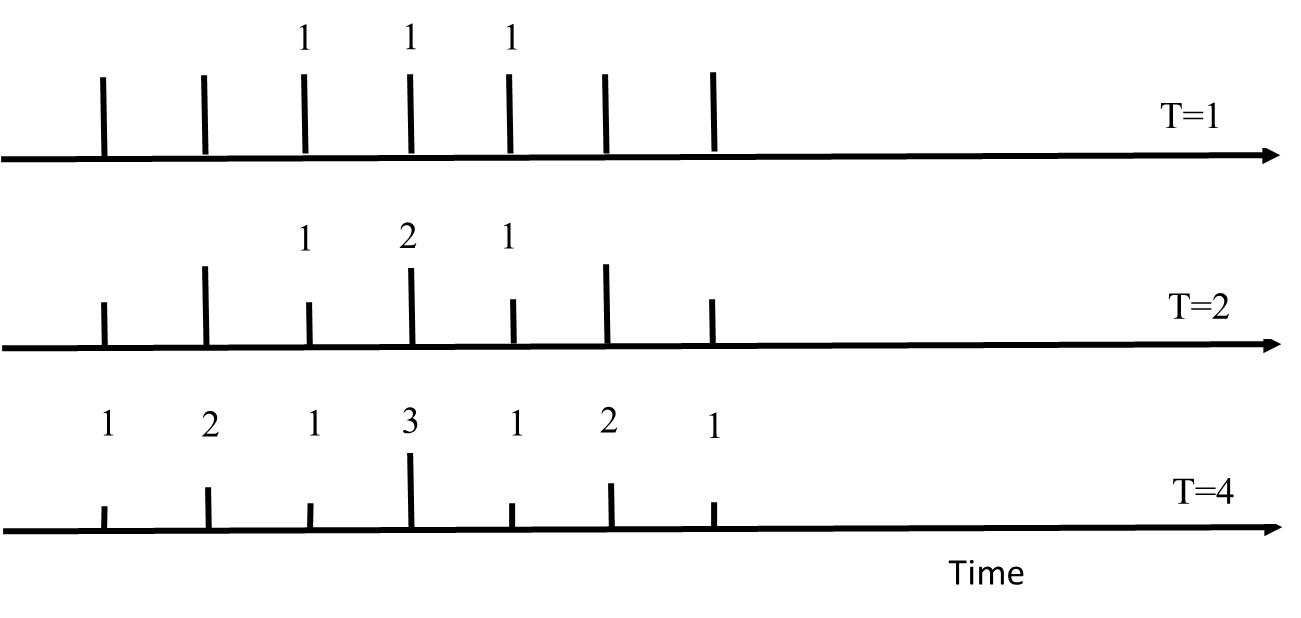}
\caption{(A) Forward horizontal visibilty of the points of a time serie of period 8 of the logistic equation. (B) As the period of the orbits doubles, the number of points in the stationary regimes also does. The forward horizonal visibility of each of these points coincides with their indices, as defined in section \ref{S1}. Thus, period doubling bifurcartions are reflected in the duplication of points of the orbits and in the value of their indices that, in this case, inform about their horizontal visibility.
 }
\label{visidesdobla}
\end{figure}

Recently, we have shown the visibility patterns behind the cascades to chaos \cite{Nuno20}.  
This paper focus on the visibility properties that appears in the bifurcation diagram of discrete maps, in particular, the unimodal maps\cite{May, Strogatz, Peitgen}.  Figure \ref{visidesdobla} depicts the horizontal visibility algorithm and the visibility pattern that is obtained as the period of the equilibrium orbits doubles. It turns out that the forward visibility sequence at the Feigenbaum accumulation point is the Gros sequence multiply by 2. Actually, the sequence coincides exactly with the sequence derived from forward horizontal visibility procedure, i.e. when the horizontal visibility toward the future of the time series is considered.

\begin{figure} 
\centering
\includegraphics[width=0.6\columnwidth]{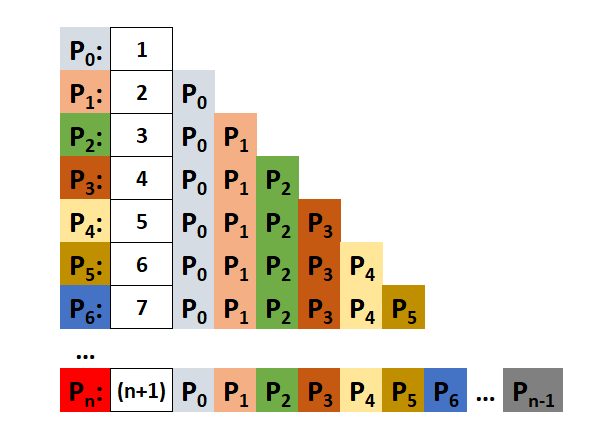}
\caption{Schematic representation of the procedure for the generation of visibility patterns in the period doubling Feigenbaum cascade by means of a recurrencce formula.  At the limit  $n \to \infty$, the visibility pattern at the onset of chaos of this cascades, $P_{\infty}$ is obtained.}
\label{Fig8}
\end{figure}

As describe in \cite{Nuno20}, the sequences $r_n$ are constructed in recursive way taking into account the period doubling of the time series. Specifically, the visibility pattern of the period $T=2^n$ can be derived from the visibility patterns of the previous periods as:  (see Figure \ref{Fig8}). 
\[
P_n = 2(n+1) \{2\} P_1 P_2 P_3 \ldots P_{n-1}  \ldots
\]
where $P_{n-1}$ is the visibility pattern for $n=1,2,\ldots$.

Notice the recurrence law that generates the $n$th-visibility pattern $P_n$ for the forward horizontal visibility: the largest visibility precedes the elementary block, $\{ 1 \}$ and after, all previous visibility patterns follow, until the one  before, $P_{n-1}$, that contains again all previous visibility patterns. 

Although it could be thought that this procedure differs from the definitions of the Gros sequence derived  in the previous section, actually, they are related, as shown in Figure \ref{visidesdobla}(B). As the period doubles, the number of points that form the orbit in the stationary regime also does. The forward visibility of each of these points is precisely their index as defined in section \ref{S1}. In other words, the point index coincides with the vibility of each of the points of the periodic orbits that occur during the Feigenbaum cascade. At the limit, at the accumulation point of this cascade, the orbit becomes chaotic and, even in the stationary regime, it is formed by an infinite number of points. Nonetheless, the visibility pattern contains the fingerprint of the period doubling cascade as it is reflected in the infinite ruler sequence.

\section{Concluding remarks}

As shown in this paper, the sequence known as the ruler function or Gros sequence appears naturally in multiple problems in fundamental mathematical fields ranging from number theory to data analysis, including classical and fractal geometry. Here, we have added four new descriptions of this sequence  to the already existing ones \cite{A001511}. All of them have in common that the ruler sequence is obtained by a recursive procedure based on the duplication of some of the points presented at each step, either intervals, vertices  or points. This is in agreement with the definition of the $n$th-term of the ruler sequene, $a(n)$, as the highest power of 2 that divides $2\,n$.   Nonetheless, it is appropiate to remark that this definition does not capture the recursive nature of the sequence and, as a consequence, its self-containing and fractal properties \cite{Kimberling}. To this respect, it should be noted that the definition of fractal sequence given in \cite{Kimberling} seems too restrictive because it is not satisfied by the ruler sequence, despite it is directly derived from the fractal structure of the Cantor set.

The dynamic description of the ruler sequence given in this paper is also interesting because, to our knowledge, it provides a first example of a demographic model that shares fractal properties. Indeed, as depicted in Fig. \ref{pyr}, the age distribution, $a$, follows the power law $p(a) = 2^{-a}$. Even more appealing is the fact that, because the self-containing property, the same age distribution is obtained if the population starts to grow from any of its individuals taken separately at a given time. We derived a recursive formula for this discrete dynamical model \ref{discreteq2} that computes the population at step $n+1$ in terms of the two previous populations $N(n)$ and $N(n-1)$. Moreover, the solution of this equation is an exponential growth $N(n) = 2^n - 1$.  As commented in section \ref{automata}, the demographic model can be generalized to consider the death of individuals at a given age. In this case, the age distribution of the population does not follow the ruler function.

It is worthy to note that the pattern that are derived in each step to construct $r_n$ is different in the Cantor set and the visibility map. The reason has to do with the definition of the process of division.   Certainly, we could also derive the Gros sequence from 
half of the Cantor set, i.e. by using only the left and middle intervals. In this case, the succesive elements of the Gros sequence would be:
\[
r_1 = \{1\}; \,   r_2 = \{1, 2\}; \, r_3 = \{1, 2, 1, 3 \}; \, r_4 = \{1, 2, 1, 3, 1, 2, 1, 4\}
\]
in the same way as it is obtained in the visibility procedure. The same can be done if we restrict the set of vertex indices to half of the sides of the polygons. Needless to say that any of the procedures yield the same final infinite sequence, the ruler sequence.

If the complete visibility algorithm is applied to obtain the visibilty patterns of the time series of the Feigenbaum cascades, the double of the ruler sequence, i.e. the sequence obtained by multiplying by 2 all the terms, is obtained \cite{Nuno20}. This is due to the temporal symmetry of the periodic series for each period of the cascades, i.e. the visibility is the same looking forward than backward. Similar horizontal vibility sequences are obtained for other periodic doubling bifurcations of unimodal maps as, for instance, the one that corresponds to the period 3 cascade \cite{A333363}.

The ruler or Gros sequence is defined in the On-Line Encyclopedia of Integer Sequences (OEIS) \cite{A001511} as an integer sequence that appears in the division of even number by powers of 2. There, it is also related to other problems, as the Tower of Hanoi and the Chinese Rings. In this paper, we have shown that the ruler sequence is ubiquitous and appears in other classical mathematical problems as: (i) population dynamics, (ii) Cantor set, (iii) polygon construction and (iv) chaotic bifurcations. The common properties behind the generation of the ruler sequence are duplication and recursiveness. We hope that this paper encourages searching for this fundamental sequence in other scientific fields that complement its characterization.


\newpage


\begin{thebibliography}{11}


\bibitem{A001511} Sloane, N. J. A. Sequence A001511 in "The On-Line Encyclopedia of Integer Sequences." 
http://www.oeis.org

\bibitem{A333363}  Sloane, N. J. A. Sequence A333363  in "The On-Line Encyclopedia of Integer Sequences." 
http://www.oeis.org


\bibitem{Allouche} Allouche, J.P., Jeffrey Shallit. {\it Automatic sequences:  theory, applications, generalizations}, Cambridge University Press (2003)



\bibitem{Barnes} Barnes, J. {\it Gems of Geometry}, Springer (2012)

\bibitem{Boyer} Boyer C.B. and U.C. Merzbach.  {\it A history of mathematics}. Wiley (2011) 

\bibitem{Caroli} Caroli, A. in: https://oeis.org/A001511 (Feb 26, 2016)

\bibitem{Caswell} Keyfitz N. and H. Caswell. {\it Applied mathematical demography}, Springer (2005)

\bibitem{Cobeli} Cobeli, C., M. Prunescu and A. Zaharescu. {\it A growth model based on the arithmetic Z -game}. Chaos, Solitons and Fractals 91, 136–147 (2016) 


\bibitem{Conway} Conway J. H. and  Richard K. Guy. {\it  The Book of Numbers}, Copernicus (1996)

\bibitem{Dunham} Dunham, W. {\it The Calculus Gallery: Masterpieces from Newton to Lebesgue},  Princeton University Press (2018)

\bibitem{Falconer} Falconer K.  {\it Fractal Geometry: Mathematical Foundations and Applications}, Wiley (2014)

\bibitem{Graham}  Graham R.L., D. E. Knuth and O. Patashnik, {\it Concrete Mathematics}. Addison-Wesley, Reading, MA, (1990)

\bibitem{Hanoi} Hinz, A.M, S. Klavzar, U. Milutinovic, C. Petr, I. Stewart. {\it The Tower of Hanoi:  Myths and Maths}, Birkhäuser (2013)

\bibitem{Kimberling} Kimberling, C. {\it Proper Self-Containing Sequences, Fractal
Sequences, and Para-Sequences}, https://oeis.org/A131987/a131987.pdf

\bibitem{Kuh}  Kuh, Devin. (2013). {\it Constructible regular $n$-gons}. Whitman College.

\bibitem{Larcombe} Larcombe, P. J., O. D. Bagdasar and E. J. Fennessey, Horadam sequences:
a survey, Bull. I.C.A., 67,  49–72. (2013)

\bibitem{Luque} Luque B, L. Lacasa, F.J. Ballesteros and A. Robledo. {\it Feigenbaum Graphs:
A Complex Network Perspective of Chaos}. PLoS ONE 6(9): e22411.
doi:10.1371/journal.pone.0022411 (2011)

\bibitem{May} May, R. M. {\it Simple mathematical models with very complicated dynamics}, Nature, 261 (1976)

\bibitem{Nuno20}  Nuño, J.C. and F. J. Muñoz.  
Universal visibility patterns of unimodal maps.  Chaos 30, 063105 (2020); https://doi.org/10.1063/5.0006652 

\bibitem{Peitgen}  Peitgen, H.O., Hartmut Jürgens, Dietmar Saupe. {\it Chaos and Fractals New Frontiers of Science}, Second Edition-Springer (2004)

\bibitem{R} R Core Team (2014). {\it R: A language and environment for statistical computing}. R Foundation for Statistical Computing, Vienna, Austria. URL http://www.R-project.org/

\bibitem{Roberts}  Roberts, J. {\it  Lure of the integers}, The Mathematical Association of America (1992)

\bibitem{Salinelli} Salinelli E., Franco Tomarelli. {\it  Discrete Dynamical Models}, Springer International Publishing (2014)

\bibitem{Strogatz} Strogatz, S.H. {\it Nonlinear Dynamics and Chaos with Applications to Physics,
Biology, Chemistry, and Engineering}, Westview Press-CRC Press (2018).


\bibitem{Vallin} Vallin, R. W.  {\it The Elements of Cantor Sets with Applications}, Wiley (2013)

\bibitem{Wolfram}  Wolfram, S. (Ed.). Theory and Application of Cellular Automata. Reading, MA: Addison-Wesley (1986)

\end{thebibliography}
\end{document}